\documentclass[11pt,reqno]{amsart}
\usepackage{amssymb}
\usepackage{amsmath}
\usepackage{amsfonts}
\usepackage{graphicx}
\usepackage[utf8]{inputenc}

\usepackage[total={17cm,23.5cm},top=2.5cm, left=2.3cm]{geometry}
\usepackage{hyperref}
    \usepackage{aeguill}
    \usepackage{type1cm}
    \usepackage[usenames,dvipsnames]{color}

\usepackage[dvipsnames]{xcolor}

\newtheorem{theorem}{Theorem}[subsection]

\newtheorem{corollary}[theorem]{Corollary}

\newtheorem{proposition}[theorem]{Proposition}
\newtheorem{conjecture}[theorem]{Conjecture}

\newtheorem{problem}[theorem]{Problem}

\newtheorem{definition}[theorem]{Definition}

\newtheorem{remark}[theorem]{Remark}

\DeclareMathOperator{\dens}{dens\,}

\DeclareMathOperator{\diam}{diam}
\DeclareMathOperator{\dist}{dist}

\newcommand{\N}{\mathbb{N}}

\newcommand{\R}{\mathbb{R}}

\newcommand{\supp}{\operatorname{supp}}

\newcommand{\espan}{\operatorname{span}}

\usepackage[dvipsnames]{xcolor}

\newcommand{\Ker}{\operatorname{Ker}}

\def\dist{{\mathop\mathrm{\,dist\,}}}

\def\bint{{\ifinner\rlap{\bf\kern.35em--}
\int\else\rlap{\bf\kern.45em--}\int\fi}\ignorespaces}

\def\bbint{{\ifinner\rlap{\bf\kern.35em--}
\hspace{0.078cm}\int\else\rlap{\bf\kern.45em--}\int\fi}\ignorespaces}

\def\diam{{\mathop\mathrm{\,diam\,}}}

\def\dfrac{\displaystyle\frac}

\def\bint{{\ifinner\rlap{\bf\kern.35em--}
\int\else\rlap{\bf\kern.45em--}\int\fi}\ignorespaces}

\begin{document}

\title[Corrigendum to ``On Smooth extensions of vector-valued functions]
{Corrigendum  to ``On Smooth extensions of vector-valued functions defined on closed subsets of Banach spaces" [Math. Ann. (2013) 355: 1201--1219]}

\author{M. Jim\'enez-Sevilla}
\address{M. Jiménez-Sevilla. \newline Instituto de Matem\'atica Interdisciplinar (IMI) and Departamento de An\'alisis Matem\'atico y Matem\'atica Aplicada, Facultad de Ciencias Matem\'aticas, Universidad Complutense de Madrid, Madrid 28040, Spain} \email{marjim@ucm.es}

%\subjclass[2020]{46B87, 46B20, 46G05, 46B03, 46E50, 46B25, 46B28.}
\date{\today}
%\keywords{}

%\thanks{{\em Research  supported } by Ministerio de Ciencia e Innovación Grant PID2022-138758NB-I00 (Spain).}

%%%%%%%%%%%%%%%%%%%%%%%%%%%%%%%%%%%%%%%%%%%%%%%%%%%%%%%%%%%%%%%%%%%%%

\begin{abstract} The present note is a corrigendum to the paper ``On Smooth extensions of vector-valued functions defined on closed subsets of Banach spaces" [Math. Ann. (2013) 355: 1201--1219].
\end{abstract}

%%%%%%%%%%%%%%%%%%%%%%%%%%%%%%%%%%%%%%%%%%%%%%%%%%%%%%%%%%%%%%%%%%%%%

\maketitle

%This note is a corrigendum of the paper \cite{JSSG}. 
L. Zajicek  kindly pointed us a flaw in the proof of Lema 3.8 in \cite{JSSG}. 
The present note is devoted to explain how to overcome the flaw in \cite{JSSG} with the additional assumption that $(X,Z)$ has property (E) (we follow the notation of \cite{JSSG}). 
We will give in this note 
a list of  changes to be made in the proof of Lemma 3.8 in \cite{JSSG} in order
to obtain a correct version of this lemma
and thus  correct versions of Theorem 3.1 and Theorem  3.2 of \cite{JSSG}. In a personal communication L. Zajicek informed us that
 together with J. Johanis they have obtained a different proof  of Theorem 3.1  and    Theorem 3.2 in \cite{JZ}  with the additional assumption that the pair of Banach spaces $(X,Z)$ has property (E) as well.
 
Let us see how to modify Lemma 3.8. In the first place, as it is already mentioned, we additionally should assume in Lemma 3.8 that
the pair of Banach spaces 
$(X,Z)$ has  property (E) for a constant $K>0$. Let us fix a constant $\tau>2$.
Following the notation of the proof of Lemma 3.8, we may assume, when selecting the open covering $\{B_\gamma:=B({y_\gamma},r_\gamma)\}_{\gamma \in \Gamma}$ satisfying all conditions specified in the proof of Lemma 3.8 that  additionally $\{B({y_\gamma},\frac{r_\gamma}{\tau})\}_{\gamma \in \Gamma}$ is an open cover of $A$ as well.
Notice that from \cite{JSSG},  $\|D(y)\|\le M$ for  all $y\in A$, $\|D(y)-h'(y)\|\le \varepsilon'$ for  all $y\in A$ and $\|h'(y)-h'(y_\gamma)\|\le \frac{\varepsilon'}{8C_0}$
for all  $y\in B_\gamma$, so we get 
\begin{equation}\label{cotaparah'}
\|h'(y)\|\le \|h'(y)-h'(y_\gamma)\|+\|h'(y_\gamma)-D(y_\gamma)\|+\|D(y_\gamma)\|
\le \frac{\varepsilon'}{8C_0}
+\varepsilon'+M, \qquad 
\text{ for all } y\in B_\gamma.
\end{equation}
Also from \cite{JSSG},
$Lip(f-h|_A)<\min\{\frac{\varepsilon'}{C_0(1+2C_0)},\, (\frac{\varepsilon'}{2C_0(R+2\varepsilon')})^2\}$
whenever $R<\infty$ and $Lip(f-h|_A)<\frac{\varepsilon'}{C_0}$
whenever $R=\infty$ (property (iii) in the proof of Lemma 3.8), so we get 
\begin{equation*}Lip(h|_A)\le \min\{\frac{\varepsilon'}{C_0(1+2C_0)},\, (\frac{\varepsilon'}{2C_0(R+2\varepsilon')})^2\}+
Lip(f):=L_1 \qquad 
\text{ whenever } R<\infty
\end{equation*}
and $
Lip(h|_A)\le \frac{\varepsilon'}{C_0}+Lip(f):=L_2 $ whenever  $R=\infty$.
Since $L_1\le L_2$, it will be enough for us to consider $L_2$  for both cases. So
$Lip(h|_A)\le L_2$ for any $R \in (0,\infty]$.

\medskip

Let us define the open balls $V_\gamma:=
B(y_\gamma,  \frac{r_\gamma}{\tau})$ for every $\gamma \in \Gamma$. We consider the open subset
$U:=\cup_{\gamma \in \Gamma}{V_\gamma}\subset X$ and the restriction mapping $h|_{U}$.
Let us check that 
$h|_U$ is Lipschitz:
If $x,y \in U$ there are indexes $\alpha$ and $\gamma$ such that $x\in V_{\alpha}$ and $y\in V_{\gamma}$. Then by inequality \eqref{cotaparah'}, we have
\begin{align*}
\|h(x)-h(y)\|\le &\|h(x)-h(y_\alpha)\|+\|h(y_\alpha)-h(y_\gamma)\|+\|h(y_\gamma)-h(y)\| \le \\ \le & (\frac{\varepsilon'}{8C_0}
+\varepsilon'+M)\|x-y_\alpha\|+L_2 \|y_\alpha-y_\gamma\|+(\frac{\varepsilon'}{8C_0}
+\varepsilon'+M)\|y-y_\gamma\|
\le \\ \le &
(\frac{\varepsilon'}{8C_0}
+\varepsilon'
+M)(\|x-y_\alpha\|+\|y-y_\gamma\|)+L_2 (\|y_\alpha-x\|+\|x-y\|+\|y-y_\gamma\|).
\end{align*}

Now, define $\delta=\frac{2}{\tau-2}>0$.
If $\|x-y_\alpha\|\le \delta\|x-y\|$ and $\|y-y_\gamma\|\le \delta\| x-y\|$, then
\begin{align*}
\|h(x)-h(y)\|\le & \bigg[2\delta
(\frac{\varepsilon'}{8C_0}
+\varepsilon'+M)+L_2(2\delta+1)\bigg]
\|x-y\|,
\end{align*}
%where $M_i=M_1$ whenever $R<\infty$ and $M_i=M_2$ whenever $R=\infty$. 
If $\|x-y_\alpha\|>\delta\|x-y\|$
or $\|y-y_\gamma\|>\delta\|x-y\|$,
let us assume for example the first inequality,
then $\|x-y\|<\frac{1}{\delta}\|x-y_\alpha\|<\frac{r_\alpha}{\delta \tau}$ and thus $\|y-y_\alpha\|\le \|y-x\|+\|x-y_\alpha\|<\frac{r_\alpha}{\delta \tau}+\frac{r_\alpha}{\tau}=(1+\frac{1}{\delta})\frac{1}{\tau}r_\alpha=\frac{r_\alpha}{2}$ and $y\in B_\alpha$. So $\|h(x)-h(y)\|\le (\frac{\varepsilon'}{8C_0}
+\varepsilon'+M)\|x-y\|$. 
A similar argument works for the second inequality.
Therefore, 
\begin{align*}
\|h(x)-h(y)\|\le L_3
\|x-y\|, \text{ for all } x,y\in U,
\end{align*}
where
\begin{equation*}
 L_3:=\max\{2\delta
(\frac{\varepsilon'}{8C_0}
+\varepsilon'+M)+L_2(2\delta+1), \, \frac{\varepsilon'}{8C_0}
+\varepsilon'+M\}.
\end{equation*}
%$M_i=M_1$ whenever $R<\infty$ and $M_i=M_2$ whenever $R=\infty$. 

The next change to be done is to define the new open cover of $X$ formed by $V_0=X\setminus A$ and $\mathcal{C}=\{V_\beta:\beta \in \Sigma=\Gamma \cup \{0\}\}$ and consider
an open refinement $\{W_{n,\beta}\}_{n\in \mathbb N, \beta \in \Sigma}$ of $\mathcal{C}:=\{V_\beta:\,\beta \in \Sigma\}$ and a $C^1$ smooth  and Lipschitz partition of unity $\{\psi_{n,\beta}:\,n\in \mathbb N, \,\beta \in \Sigma\}$ subordinated to $\mathcal{C}$ satisfying properties (P1)-(P3) specified in the proof of Lemma 3.7 in \cite{JSSG}, where among other properties, $\supp \psi_{n,\beta}\subset W_{n,\beta}\subset V_\beta$,
for all $n\in \mathbb N$ and $\beta \in \Sigma$.

The definitions of the constants $L_{n,\beta}:=\max\{Lip(\psi_{n,\beta}),1\}$ for
$n\in \mathbb N$ and $\beta\in \Sigma$ and the functions $\delta_{n,\gamma}$ obtained by applying property (*) to $T_\gamma-h$ on $B_\gamma$ for every $n\in \mathbb N$, $\gamma \in \Gamma$ remain the same as in \cite{JSSG}.

The next modification is  to apply property (E)
to $h|_U:U\rightarrow B_{Z}(0,R+\varepsilon'+\frac{\varepsilon'}{8C_0})$.
In the case, $R=+\infty$ we get
a $KL_3$-Lipschitz function $H:X\rightarrow Z$ such that
$H|_U=h|_U$. We apply property (*) to $H$ in order to obtain
$C^1$ smooth mappings
$F_0^n:X\rightarrow Z$
such that $\|H(x)-F^n_0(x)\|\le \frac{\varepsilon'}{2^{n+2}L_{n,0}}$ for all $x\in X$ and
$
Lip(F^n_0)\le C_0KL_3$ for all $n\in \mathbb N$.
In the case $R<+\infty$ we consider the auxiliary function
$G:U\cup W\rightarrow B_{Z}(0,R+\varepsilon'+\frac{\varepsilon'}{8C_0})$, where 
$W=\{x\in X:\dist(x,U)\ge
(R+\varepsilon'+\frac{\varepsilon'}{8C_0})L_3^{-1}\}$
defined by $G|_U=h|_U$ and
$G|_W=0$. It can be checked that $G$ is $L_3$-Lispchitz. By property (E)
there is a $KL_3$-Lispchitz extension $L:X\rightarrow Z$ of
$G$ and thus, in particular, it can be checked that
$L:X\rightarrow B_Z(0, (K+1)(R+\varepsilon'+\frac{\varepsilon'}{8C_0}))$.
Now, let us apply property (*) to $L$ to get
$C^1$ smooth mappings $F_0^n:X\rightarrow B_Z(0, (K+1)(R+\varepsilon'+\frac{\varepsilon'}{8C_0})+\frac{\varepsilon'}{2^{n+2}L_{n,0}})$ such that $\|L(z)-F^n_0(z)\|<\frac{\varepsilon'}{2^{n+2}L_{n,0}}$ for all 
$z\in X$
and $
Lip(F^n_0)\le C_0KL_3$.

The definitions of the functions $\{\Delta^n_\beta\}_{n\in \mathbb N, \beta \in \Sigma }$  in terms of the functions $\{F^n_0\}_{n\in \mathbb N}$, $\{T_\gamma-\delta_{n,\gamma}\}_{n\in \mathbb N, \gamma\in \Gamma}$ 
and the definition of $g$ in terms of 
$\{\psi_{n,\beta}\}_{n\in \mathbb N, \beta \in \Sigma}$ is the same as in   Lemma 3.8 in  \cite{JSSG}.
The new upper bound for $g$ is
\begin{equation} \label{cotaparag} \tag{C.1}
    \|g(x)\|\le \sum_{(n,\beta)\in \mathbb N\times\Sigma}
    \psi_{n,\beta}(x)\|
\Delta_{\beta}^n(x)\|\le
(K+1)(R+\varepsilon'+\frac{\varepsilon'}{8C_0})+\frac{\varepsilon'}{8}:=Q_1, \quad \text{ for }x\in X.
\end{equation}
It can be checked that  
the new upper bound for
$\|(\Delta^n_\beta)'(x)\|$ for $x\in X$    is
\begin{equation*}
\|(\Delta^n_\beta)'(x)\|
\le  \max\{Lip(F^n_0), \,M+9\frac{\varepsilon'}{8}\}\le \max\{ C_0K L_3, \,
 M+9\frac{\varepsilon'}{8}\} =C_0KL_3,
\end{equation*}
where $n\in \mathbb N$ and $\beta \in \Sigma$.
%where $M_i=M_1$ for $R<+\infty$ and $M_i=M_2$ for $R=+\infty$.
Let us check that $g$
is Lispchitz. Let us denote $S:=H$ if $R=+\infty$ and $S:=L$ if $R<+\infty$.
Since $\supp \psi_{n,\beta}\subset V_\beta$ for all $n\in \mathbb N$, $\beta \in \Gamma$ and $S|_U=h$,  
following the notation of the proof of Lemma 3.8 in  \cite{JSSG} we have for $x\in X,$
\begin{align} \notag
&||g'(x)||\le \sum_{(n,\beta)\in F_x}
\|\psi'_{n,\beta}(x)\|\|S(x)-\Delta_\beta^n(x)\|+
\sum_{(n,\beta)\in F_x}\psi_{n,\beta}(x)\|(\Delta^n_\beta)'(x)\|\le \\
&\le  \sum_{\{n:(n,\beta(n))\in F_x\}} L_{n,\beta(n)}\,\frac{\varepsilon'}{2^{n+2}L_{n,\beta(n)}}+\sum_{\{n:(n,\beta(n))\in F_x\}}\psi_{n,\beta}(x)C_0KL_3\le \frac{1}{4}\varepsilon'+C_0KL_3:=L_4, \tag{C.2}
\label{cotaparag'}
\end{align}
%for $i=1,2$ where $P_i=P_1$ whenever $R<+\infty$ and $P_i=P_2$ whenever $R=\infty$.
In addition, properties (i), (ii) and (iii) in the statement of Lemma 3.8 in \cite{JSSG} follow in the same way.

The new upper bound in property (iv) in the statement of Lemma 3.8 in \cite{JSSG} is $Q_1$ given in \eqref{cotaparag}. %Notice that $Q_1$ is a function $Q_1(\varepsilon',R)$
%and $Q_1(\varepsilon',R)\to 0$ whenever $\varepsilon'\to 0, R\to 0$.
The new upper bound in property (v) in the statement of Lemma 3.8 in \cite{JSSG} is $L_4$
%$R_1$ for $R<+\infty$ and $R_2$ for $R=+\infty$ 
given in \eqref{cotaparag'}.
%Notice that $R_1$ is a function $R_1(\varepsilon',R,M, Lip(f))$
%and $R_1(\varepsilon',R,M, Lip(f))\to 0$ whenever $\varepsilon'\to 0, R\to 0, M\to 0, Lip(f)\to 0$. Also,  $R_2$ is a function $R_2(\varepsilon',M, Lip(f))$
%and $R_2(\varepsilon',M, Lip(f))\to 0$ whenever $\varepsilon'\to 0,M\to 0, Lip(f)\to 0$.

Due to the preceding arguments, we must assume the additional condition that the pair of Banach spaces $(X,Z)$ has property (E) in the statements of Theorem 3.1 and Theorem 3.2 of \cite{JSSG}.
In order to apply Lemma 3.8  with the new upper bounds $Q_1$ and $L_4$, 
let us consider simpler upper bounds in terms of $M,K,C_0,\,\varepsilon$.
Recall that, from the assumptions in Lemma 3.8 in \cite{JSSG}, $3\varepsilon'<\varepsilon$,
and thus after some straightforward calculations it can be checked that 
\begin{align*}
Q_1&< (K+1)R+A_1\varepsilon, \\
L_4&< A_2\varepsilon+C_0K
(4\delta+1)\max\{Lip(f),M\},
\end{align*}
where $A_1:=\frac{1}{3}[(K+1)(1+\frac{1}{8C_0})+\frac{1}{8}]$,\,
 $A_2:=\frac{1}{3}[\frac{1}{4}+C_0K(1+\frac{2}{C_0}+\frac{1}{8C_0})]$ and
we have assumed without loss of generality that
 $2\delta\le 1$.

\medskip

Now, we can reproduce the proofs of 
Theorem 3.1 and 3.2  given in \cite{JSSG} 
 with any sequence
$\{\varepsilon_n\}_n\subset (0,1)$ such that $\sum_n \varepsilon_n <\infty$.
For the sake of completeness let us briefly sketch it. Following the notation of the proofs of Theorem 3.1 and 3.2 in \cite{JSSG}, 
in the first step we apply Lemma 3.7   to get a $C^1$
smooth mapping $G_1:X\rightarrow Z$ such that if $g_1:=G_1|_A$ then
(i) $||f(y)-g_1(y)||<\varepsilon_1$ for all $y\in A$; (ii) $\|D(y)-G_1'(y)\|<\varepsilon_1$
for all $y\in A$ and 
(iii) $Lip(f-g_1)<\varepsilon_1$.
In the second step, the  funtion $f-g_1:A\rightarrow B_Z(0,\varepsilon_1)$ satisfies the mean value condition for the function $D-G_1':A\rightarrow \mathcal{L}(X,Z)$ with
$\|D(y)-G_1'(y)\|\le \varepsilon_1$ for all
$y\in A$. So Lemma 3.8
applies to get a $C^1$ smooth function $G_2:X\rightarrow Z$ such that if $g_2:=G_2|_A$,
then (i) $\|f(y)-g_1(y)-g_2(y)\|<\varepsilon_2$
 for all $y\in A$,
(ii)  $\|D(y)-G_1'(y)-G_2'(y)\|<\varepsilon_2$
for all $y\in A$, 
(iii) $Lip(f-g_1-g_2)<\varepsilon_2$,
(iv) $\|G_2(y)\|<(K+1)\varepsilon_1+A_1\varepsilon_2$ for all
$y\in X$ and (iv) $Lip(G_2)< A_2\varepsilon_2+C_0K(4\delta+1)\varepsilon_1 $.
In general, in the $n$-th step ($n\ge 2$) we apply Lemma 3.8 to the Lipschitz function $f-(g_1+\cdots+g_{n-1}):A\rightarrow B_Z(0,\varepsilon_{n-1})$, which satisfies
 the mean value condition for the function $D-(G_1'+\cdots+G_{n-1}'):A\rightarrow \mathcal{L}(X,Z)$ with $\|D(y)-(G_1'+\cdots+G_{n-1}')(y)\|\le \varepsilon_{n-1}$ for $y\in A$ and $Lip (f-(g_1+\cdots+g_{n-1}))<\varepsilon_{n-1}$ to get a $C^1$ smooth function $G_n:X\rightarrow Z$ such that if $g_n:=G_n|_A$,
then (i) $\|f(y)-(g_1+\cdots+g_n)(y)\|<\varepsilon_n$ for all $y\in A$,
(ii)  $\|D(y)-(G_1'+\cdots +G_n')(y)\|<\varepsilon_n$
for all $y\in A$, 
(iii) $Lip(f-(g_1+\cdots+g_n))<\varepsilon_n$,
(iv) $\|G_n(y)\|<(K+1)\varepsilon_{n-1}+A_1\varepsilon_n$ for 
all $y\in X$, and  (v) $Lip(G_n)< A_2\varepsilon_n+C_0K(1+4\delta)
\varepsilon_{n-1}$.

\medskip

Then, it can be checked that the sum $G:=\sum_nG_n$ is a $C^1$ smooth extension of $f$.
If the initial function $f$ is Lipschitz and satisfies the mean value condition for a function $D$ with $\|D(y)\|\le M$ for all $y\in A$, then applying Lemma 3.8 instead of Lemma 3.7 in the first step, we additionally obtain  $Lip(G_1)< A_2\varepsilon_1+C_0K(1+4\delta)\max\{M, Lip(f)\}$. So $$Lip(G)< C_0K(1+4\delta)\max\{M,Lip(f)\}+
\widetilde{\varepsilon},$$
where $\widetilde{\varepsilon}=A_2\varepsilon_1+\sum_{n\ge 2}(A_2\varepsilon_n+C_0K(1+4\delta)
\varepsilon_{n-1})$. Since  we can choose from the beginning $\widetilde{\varepsilon}>0$
and $\delta>0$ to be as small as we want, we can get $G$ in such a way that  $Lip(G)$ is as close to $C_0K\max\{M,Lip(f)\}$
as we need.

So in the statement of Theorem 3.2 in \cite{JSSG} the upper bound for $Lip(G)$ should be replaced by any constant greater than  $C_0K\max\{M,Lip(f)\}$.

\medskip

In addition, the examples given in Corollary 3.4 in \cite{JSSG} are valid since they have property (E) and (*) (or equivalently property (E) and (A)). However,  the upper bound for $Lip(G)$ should be replaced by the new one obtained above.

\medskip

Finally, let us metion that  A. Sofi kindly pointed us that Corollary 4.11 does not hold. Although  (i)$\Rightarrow$(ii)$\Rightarrow$(iii)$\Rightarrow$ $X^{**}$ is a $\mathcal{P}_{\lambda}$-space, the implication  (iv)$\Rightarrow$(i)  does not hold.

%Also consider the auxiliary function $F:U\cup W\rightarrow Z$ with
%$F|_U=1$ and $F|_W=0$
%Since $X$ has property (*) there is 

\medskip

{\bf Conflict of interest}. The author states that there is no conflict of interest.

{\bf Data availability.}
This note has no additional data.

\end{document}